\documentclass[11pt,a4paper]{article}
\usepackage{mathrsfs}
\usepackage{amsmath}
\usepackage{amsfonts}
\usepackage{amssymb}

\usepackage{graphicx,subfigure,epsfig}
\usepackage{enumerate}
\usepackage{arydshln}
\usepackage[pdftex,linkcolor=blue,citecolor=blue,backref=pagebackref]{hyperref}
\usepackage{color}
\usepackage{pdfpages}
\usepackage{subfigure}
\usepackage{tabularx}
\allowdisplaybreaks[4]

\begin{document}

\renewcommand{\thefootnote}{\fnsymbol{footnote}}

\newtheorem{theorem}{Theorem}[section]
\newtheorem{corollary}[theorem]{Corollary}
\newtheorem{definition}[theorem]{Definition}
\newtheorem{conjecture}[theorem]{Conjecture}
\newtheorem{question}[theorem]{Question}
\newtheorem{lemma}[theorem]{Lemma}
\newtheorem{proposition}[theorem]{Proposition}
\newtheorem{example}[theorem]{Example}
\newtheorem{fact}[theorem]{Fact}
\newenvironment{proof}{\noindent {\bf
Proof.}}{\rule{3mm}{3mm}\par\medskip}
\newcommand{\remark}{\medskip\par\noindent {\bf Remark.~~}}
\newcommand{\pp}{{\it p.}}
\newcommand{\de}{\em}
\newcommand{\qbinom}[2]{\genfrac{[}{]}{0pt}{}{#1}{#2}}

\newcommand{\JEC}{{\it Europ. J. Combinatorics},  }
\newcommand{\JCTB}{{\it J. Combin. Theory Ser. B.}, }
\newcommand{\JCT}{{\it J. Combin. Theory}, }
\newcommand{\JGT}{{\it J. Graph Theory}, }
\newcommand{\ComHung}{{\it Combinatorica}, }
\newcommand{\DM}{{\it Discrete Math.}, }
\newcommand{\ARS}{{\it Ars Combin.}, }
\newcommand{\SIAMDM}{{\it SIAM J. Discrete Math.}, }
\newcommand{\SIAMADM}{{\it SIAM J. Algebraic Discrete Methods}, }
\newcommand{\SIAMC}{{\it SIAM J. Comput.}, }
\newcommand{\ConAMS}{{\it Contemp. Math. AMS}, }
\newcommand{\TransAMS}{{\it Trans. Amer. Math. Soc.}, }
\newcommand{\AnDM}{{\it Ann. Discrete Math.}, }
\newcommand{\NBS}{{\it J. Res. Nat. Bur. Standards} {\rm B}, }
\newcommand{\ConNum}{{\it Congr. Numer.}, }
\newcommand{\CJM}{{\it Canad. J. Math.}, }
\newcommand{\JLMS}{{\it J. London Math. Soc.}, }
\newcommand{\PLMS}{{\it Proc. London Math. Soc.}, }
\newcommand{\PAMS}{{\it Proc. Amer. Math. Soc.}, }
\newcommand{\JCMCC}{{\it J. Combin. Math. Combin. Comput.}, }
\newcommand{\GC}{{\it Graphs Combin.}, }
\newcommand{\LAA}{{\it Linear Algeb. Appli.}, }

\title{ \bf Snevily's Conjecture about $\mathcal{L}$-intersecting Families
\\on Set Systems and its Analogue on Vector Spaces}

\author{Jiuqiang Liu$^{a,b,} \thanks{The corresponding author}$, Guihai Yu$^{a, *}$, Lihua Feng$^{c}$, and Yongjiang Wu$^{c}$\\
{\small  $^a$ College of Big Data Statistics, Guizhou University of Finance and Economics}\\
 {\small Guiyang, Guizhou, 550025, China}\\
 {\small $^{b}$ Department of Mathematics, Eastern Michigan University}\\
{\small Ypsilanti, MI 48197, USA}\\
{\small $^{c}$ School of Mathematics and Statistics, HNP-LAMA, Central South University}\\
{\small Changsha, Hunan, 410083, China}\\
{\small E-mail: { \tt jiuqiang68@126.com, yuguihai@mail.gufe.edu.cn, fenglh@163.com}}}

\maketitle
\vspace{-0.5cm}

\begin{abstract}
The classical Erd\H{o}s-Ko-Rado theorem on the size of an intersecting family of $k$-subsets of the set $[n] = \{1, 2, \dots, n\}$ is one of the fundamental intersection theorems for set systems. After the establishment of the EKR theorem, many intersection theorems on set systems have appeared in the literature, such as the well-known Frankl-Wilson theorem, Alon-Babai-Suzuki theorem, and Grolmusz-Sudakov theorem. In 1995, Snevily proposed the conjecture that the upper bound for the size of an $\mathcal{L}$-intersecting family of subsets of $[n]$ is ${{n} \choose {s}}$ under the condition $\max \{l_{i}\} < \min \{k_{j}\}$, where $\mathcal{L} = \{l_{1}, \dots, l_{s}\}$ with $0 \leq l_{1} < \cdots < l_{s}$ and $k_{j}$ are subset sizes in the family. In this paper, we prove that Snevily's conjecture holds for $n \geq {{k^{2}} \choose {l_{1}+1}}s + l_{1}$, where $k$ is the maximum subset size in the family. We then derive an analogous result for $\mathcal{L}$-intersecting families of subspaces of an $n$-dimensional vector space over a finite field $\mathbb{F}_{q}$.
\end{abstract}

 \vskip 0.1cm

\section{Introduction}
\hspace*{0.5cm}
Throughout the paper, we use $X$ for the set  $[n]=\{1, 2, \dots,
n\}$. A family $\mathcal{F}$ of subsets of $[n]$ is called
$intersecting$ if every pair of distinct subsets $E, F \in
\mathcal{F}$ have a nonempty intersection. Let $\mathcal{L}=\{l_1,
l_2, \dots, l_s\}$ be a set of $s$ nonnegative integers with $0 \leq l_{1} < \cdots < l_{s}$. A family
$\mathcal{F}$ of subsets of $[n]$ is called
$\mathcal{L}$-$intersecting$ if $|E \cap F| \in \mathcal{L}$ for
every pair of subsets $E$, $F$ in $\mathcal{F}$.
$\mathcal{F}$ is $k$-{\em uniform} if it is a collection of $k$-subsets
of $[n]$. Thus, a $k$-uniform intersecting family is
$\mathcal{L}$-intersecting for $\mathcal{L}=\{1, 2, \dots, k-1\}$.

In 1961, Erd\H{o}s, Ko, and Rado \cite{ekr} proved the following classical result.

\begin{theorem}\label{thm1.1} (Erd\H{o}s, Ko, and Rado, \cite{ekr}).
Let $n \geq 2k$ and let $\mathcal{A}$ be a $k$-uniform intersecting
family of subsets of $[n]$. Then $|\mathcal{A}|\leq {{n-1} \choose
{k-1}}$ with equality only when $\mathcal{F}$ consists of all
$k$-subsets containing a common element.
\end{theorem}

Since then, many intersection theorems have appeared in the literature.
Recall that a $k$-uniform intersecting family is
$\mathcal{L}$-intersecting for $\mathcal{L}=\{1, 2, \dots, k-1\}$.
In 1975, D. R. Ray-Chaudhuri and R. M. Wilson derived the following
tight upper bound for a $k$-uniform $\mathcal{L}$-intersecting family, where $\mathcal{L}$ is an arbitrary set of $s$ non-negative integers.

\begin{theorem}\label{thm1.2} (Ray-Chaudhuri and Wilson, \cite{rw}).
Let $\mathcal{L}=\{l_1, l_2, \dots, l_s\}$ be a set of $s$
nonnegative integers. If $\mathcal{A}$ is a $k$-uniform
$\mathcal{L}$-intersecting family of subsets of $[n]$, then
$|\mathcal{A}|\leq {{n} \choose {s}}$.
\end{theorem}

In terms of parameters $n$ and $s$, this inequality is best possible, as shown by the set of all $s$-subsets of $[n]$ with
$\mathcal{L}=\{0, 1, \dots, s - 1\}$.

In 1981, Frankl and Wilson \cite{fw} provided the following celebrated theorem which extends Theorem \ref{thm1.2} by allowing different subset sizes.

\begin{theorem}\label{thm1.3} (Frankl and Wilson, \cite{fw}). Let
$\mathcal{L}=\{l_1, l_2, \dots, l_s\}$ be a set of $s$ nonnegative
integers. If $\mathcal{A}$ is an $\mathcal{L}$-intersecting family
of subsets of $[n]$, then
\[|\mathcal{A}|\leq {{n} \choose {s}}+ {{n}
\choose {s-1}}+ \cdots + {{n} \choose {0}}.\]
\end{theorem}

The upper bound in Theorem \ref{thm1.3} is best possible, as demonstrated by the set of all subsets of size at most $s$ of $[n]$. The original proof of Theorem 1.3 by Frankl and Wilson \cite{fw} uses the method of higher incidence matrices. The proof given by Alon, Babai, and Suzuki \cite{abs} is a very elegant linear algebra method which has been modified to produce many other interesting intersection theorems.
Alon, Babai, and Suzuki \cite{abs} also used this linear algebra method to derive the following theorem
which is a common generalization to Theorems \ref{thm1.2} and \ref{thm1.3}, under the convention that ${{a} \choose {b}} = 0$ if $b < 0$.

\begin{theorem}\label{thm1.4} (Alon, Babai, and Suzuki, \cite{abs}). Let
$\mathcal{L}=\{l_1, l_2, \dots, l_s\}$ be a set of $s$ nonnegative
integers and $K = \{k_1, k_2, \dots, k_r \}$ be a set of integers
satisfying $k_i > s-r$ for every $i$.
Suppose that $\mathcal{A} = \{A_1, A_2, \dots, A_{m}\}$ is a family of subsets of $[n]$ such that $|A_i| \in K$ for every $1 \leq i \leq m$
and $|A_{i} \cap A_{j}| \in \mathcal{L}$ for every pair $i \neq j$.
Then
\[m \leq {{n} \choose {s}}+ {{n} \choose {s-1}}+ \cdots + {{n} \choose {s-r+1}}.\]
\end{theorem}

This result is also best possible, as shown by the set of all subsets of $[n]$ of sizes at least $s-r+1$ and at most $s$.
In 2002, Grolmusz and Sudakov \cite{gs} proved the next theorem which extends Theorems \ref{thm1.3} and \ref{thm1.4}
to $t$-wise $\mathcal{L}$-intersecting families for $t \geq 2$.

\begin{theorem}\label{thm1.5} (Grolmusz and Sudakov, \cite{gs}). Let $t \geq 2$ and let
$\mathcal{L}=\{l_1, l_2, \dots, l_s\}$ be a set of $s$ nonnegative
integers. If $\mathcal{A} = \{A_1, A_2, \dots, A_{m}\}$ is a family
of subsets of $[n]$ such that $|A_{i_1} \cap A_{i_2}\cap \dots \cap A_{i_t}| \in \mathcal{L}$ for every collection of $t$ distinct subsets in $\mathcal{A}$. Then
\[m \leq (t - 1)\left[{{n}\choose {s}}+ {{n} \choose {s-1}}+
\cdots + {{n} \choose {0}}\right].\]
If in addition the size of every member of $\mathcal{A}$ belongs to the set $K = \{k_1, k_2, \dots, k_r\}$ and $k_{i} > s - r$ for every $i$, then
\[m \leq (t - 1)\left[{{n}\choose {s}}+ {{n} \choose {s-1}}+
\cdots + {{n} \choose {s-r+1}}\right].\]
\end{theorem}

\vspace{3mm}

Note that the set $\mathcal{L}$ in the above theorems may contain $0$. When $\mathcal{L}$ is restricted to be a set of positive integers, Snevily \cite{s1} proved the following result which implies Frankl-Wilson theorem (Theorem \ref{thm1.3}).

\begin{theorem}\label{thm1.6} (Snevily, \cite{s1}). Let
$\mathcal{L}=\{l_1, l_2, \dots, l_s\}$ be a set of $s$ positive
integers. If $\mathcal{A}$ is an $\mathcal{L}$-intersecting family
of subsets of $[n]$, then
\[|\mathcal{A}|\leq {{n-1} \choose {s}}+ {{n-1}
\choose {s-1}}+ \cdots + {{n-1} \choose {0}}.\]
\end{theorem}

In 1995, Snevily \cite{s} proposed the next conjecture which would improve the bounds in Theorems \ref{thm1.3}, \ref{thm1.4}, and \ref{thm1.6} considerably under the condition $\max \{l_{i}\} < \min \{k_{j}\}$.

\begin{conjecture}\label{conj1.7} (Snevily, \cite{s}). Let
$\mathcal{L}=\{l_1, l_2, \dots, l_s\}$ be a set of $s$ nonnegative
integers and $K = \{k_1, k_2, \dots, k_r\}$ be a set of $r$ positive integers with $\max \{l_{i}\} < \min \{k_{j}\}$. If $\mathcal{A} = \{A_1, A_2, \dots, A_{m}\}$ is an $\mathcal{L}$-intersecting family of subsets of $[n]$ such that $|A_{i}| \in K$ for every $1 \leq i \leq m$, then
\[|\mathcal{A}|\leq {{n} \choose {s}}.\]
\end{conjecture}

Snevily \cite{s} proved that Conjecture \ref{conj1.7} holds when $\mathcal{L}=\{0, 1, \dots, s - 1\}$ as shown in the following theorem and Gao et al. \cite{glx} proved that Conjecture \ref{conj1.7} is true when $n$ is sufficiently large.

\begin{theorem}\label{thm1.8} (Snevily, \cite{s}).
Let $\mathcal{L}=\{0, 1, \dots, s - 1\}$ and $K = \{k_1, k_2, \dots, k_r\}$ be a set of $r$ positive integers with $s \leq \min \{k_{j}\}$. Suppose that $\mathcal{A} = \{A_1, A_2, \dots, A_{m}\}$ is an $\mathcal{L}$-intersecting family of subsets of $[n]$ such that $|A_{i}| \in K$ for every $1 \leq i \leq m$.
Then
\[|\mathcal{A}|\leq {{n} \choose {s}}.\]
\end{theorem}

Here, we prove the next theorem which shows that Conjecture \ref{conj1.7} holds when $n \geq {{k^{2}} \choose {l_{1}+1}}s + l_{1}$.

\begin{theorem}\label{thm1.9}
Let $\mathcal{L} = \{l_1, l_2, \dots, l_s\}$ be a set of $s$ nonnegative integers with $0 \leq l_{1} < l_{2} < \cdots < l_{s}$. Suppose that $\mathcal{A} = \{A_1, A_2, \dots, A_{m}\}$ is an $\mathcal{L}$-intersecting family of subsets of $[n]$
such that $|A_{i}| \not \in \mathcal{L}$ for every $1 \leq i \leq m$.
If $n \geq {{k^{2}} \choose {l_{1}+1}}s + l_{1}$, where $k = \max_{1 \leq j \leq m} |A_{j}|$, then
\[|\mathcal{A}|\leq {{n} \choose {s}}.\]
\end{theorem}

The following result extends Theorem \ref{thm1.9} to $t$-wise $\mathcal{L}$-intersecting families which improves the bounds in Theorem \ref{thm1.5} considerably.

\begin{theorem}\label{thm1.10}
Let $\mathcal{L} = \{l_1, l_2, \dots, l_s\}$ be a set of $s$ nonnegative integers with $0 \leq l_{1} < l_{2} < \cdots < l_{s}$. Suppose that $\mathcal{A} = \{A_1, A_2, \dots, A_{m}\}$ is a family
of subsets of $[n]$ such that $|A_{i_1} \cap A_{i_2}\cap \dots \cap A_{i_t}| \in \mathcal{L}$ for every collection of $t$ distinct subsets in $\mathcal{A}$
and $|A_{i}| \not \in \mathcal{L}$ for every $1 \leq i \leq m$.
If $n \geq {{k^{2}} \choose {l_{1}+1}}s + l_{1}$,  where $k = \max_{1 \leq j \leq m} |A_{j}|$, then
\[|\mathcal{A}|\leq (t - 1){{n} \choose {s}}.\]
\end{theorem}

Next, denote
\begin{equation}
\qbinom{n}{k} = \prod_{0 \leq i \leq k - 1}\frac{q^{n - i} - 1}{q^{k - i} - 1}.
\label{1.1}\end{equation}
As for vector spaces over a finite field $\mathbb{F}_{q}$, it is well known that the number of all $k$-dimensional
subspaces of an $n$-dimensional vector space over $\mathbb{F}_{q}$ is equal to
$\qbinom{n}{k}$.

There are results on the sizes of families of subspaces of an $n$-dimensional vector space analogous to theorems on $\mathcal{L}$-intersecting families of subsets of $[n]$. The next theorem is proved by  Deza and Frankl \cite{df}, Greene and Kleitman \cite{gk}, and Hsieh \cite{h}.

\begin{theorem}\label{thm1.11} (Deza and Frankl \cite{df}, Greene and Kleitman \cite{gk}, Hsieh \cite{h}).
Let $n \geq 2k$ and $\mathbb{F}_{q}$ be a finite field of order $q$. Suppose that $\mathcal{V}$ is a collection of $k$-dimensional subspaces of an $n$-dimensional
vector space over $\mathbb{F}_{q}$ satisfying that $dim(V_{i} \cap V_{j}) > 0$ for any distinct subspaces $V_{i}$ and $V_{j}$ in $\mathcal{V}$.
Then
\[|\mathcal{V}| \leq \qbinom{n-1}{k-1}.\]
\end{theorem}

Frankl and Graham \cite{fg} derived the following theorem on $\mathcal{L}$-intersecting family of $k$-dimensional subspaces.

\begin{theorem}\label{thm1.12} (Frankl and Graham \cite{fg}).
Let $\mathcal{L}=\{l_1, l_2, \dots, l_s\}$ be a set of $s$ nonnegative integers and $\mathbb{F}_{q}$ be a finite field of order $q$. Suppose that $\mathcal{V}$ is a collection of $k$-dimensional subspaces of an $n$-dimensional
vector space over $\mathbb{F}_{q}$ satisfying that $dim(V_{i} \cap V_{j}) \in \mathcal{L}$ for any distinct subspaces $V_{i}$ and $V_{j}$ in $\mathcal{V}$.
Then
\[|\mathcal{V}| \leq \qbinom{n}{s}.\]
\end{theorem}

In 1990, Lefmann \cite{l} proved the following $\mathcal{L}$-intersecting theorem for ranked finite lattices,
and in 2001, Qian and Ray-Chaudhuri \cite{qr} extended to quasi-polynomial semi-lattices.

\begin{theorem}\label{thm1.13} (Lefmann \cite{l}, Qian and Ray-Chaudhuri \cite{qr}).
Let $\mathcal{L}=\{l_1, l_2, \dots, l_s\}$ be a set of $s$ nonnegative integers and $\mathbb{F}_{q}$ be a finite field of order $q$. Suppose that $\mathcal{V}$ is a collection of subspaces of an $n$-dimensional vector space over $\mathbb{F}_{q}$ satisfying that $dim(V_{i} \cap V_{j}) \in \mathcal{L}$ for any distinct subspaces $V_{i}$ and $V_{j}$ in $\mathcal{V}$.
Then
\[|\mathcal{V}| \leq
\qbinom{n}{s}
+ \qbinom{n}{s-1} + \cdots + \qbinom{n}{0}.\]
\end{theorem}

Alon et al. \cite{abs} derived the next result.

\begin{theorem}\label{thm1.14} (Alon, Babai, and Suzuki \cite{abs}).
Let $\mathcal{L}=\{l_1, l_2, \dots, l_s\}$ be a set of $s$ nonnegative integers and $\mathbb{F}_{q}$ be a finite field of order $q$. Suppose that $\mathcal{V}$ is a collection of subspaces of an $n$-dimensional vector space over $\mathbb{F}_{q}$ satisfying that $dim(V_{i} \cap V_{j}) \in \mathcal{L}$ for any distinct subspaces $V_{i}$ and $V_{j}$ in $\mathcal{V}$ and $dim(V_{i}) \in \{k_{1}, k_{2}, \dots, k_{r}\}$ with $k_{i} > s - r$ for every $i$.
Then
\[|\mathcal{V}| \leq
\qbinom{n}{s} + \qbinom{n}{s-1} + \cdots + \qbinom{n}{s-r+1}.\]
\end{theorem}

We say that a family $\mathcal{V}$ of subspaces is {\em Sperner} (or antichain) if no subspace is contained in another subspace in $\mathcal{V}$.
We will derive the following two theorems which strengthen Theorems  \ref{thm1.12} -- \ref{thm1.14} considerably.

\begin{theorem}\label{thm1.15}
Let $\mathcal{L}=\{l_1, l_2, \dots, l_s\}$ be a set of $s$ positive integers with $0 < l_{1} < l_{2} < \cdots < l_{s}$ and $\mathbb{F}_{q}$ be a finite field of order $q$. Suppose that $\mathcal{V}$ is a collection of subspaces of an $n$-dimensional vector space over $\mathbb{F}_{q}$ satisfying that $dim(V_{i} \cap V_{j}) \in \mathcal{L}$ for any distinct subspaces $V_{i}$ and $V_{j}$ in $\mathcal{V}$.
If $n \geq \max\bigg\{\log_{q}\left((q^{s}-1)
\left[ k^{2} \atop l_{1} + 1   \right] + 1 \right)  + l_{1}, \mbox{ } 2s + 1\bigg\}$, where $k = \max_{V \in \mathcal{V}} dim(V)$,
then
\[|\mathcal{V}| < \qbinom{n}{s}.\]
\end{theorem}

\begin{theorem}\label{thm1.16}
Let $\mathcal{L}=\{l_1, l_2, \dots, l_s\}$ be a set of $s$ nonnegative integers with $0 \leq l_{1} < l_{2} < \cdots < l_{s}$ and $\mathbb{F}_{q}$ be a finite field of order $q$. Suppose that $\mathcal{V}$ is a Sperner family of subspaces of an $n$-dimensional vector space over $\mathbb{F}_{q}$ satisfying that $dim(V_{i} \cap V_{j}) \in \mathcal{L}$ for any distinct subspaces $V_{i}$ and $V_{j}$ in $\mathcal{V}$.
If $n \geq \max\bigg\{\log_{q}\left((q^{s}-1)
\left[ k^{2} \atop l_{1} + 1   \right] + 1 \right)  + l_{1}, \mbox{ } 2s + 1\bigg\}$, where $k = \max_{V \in \mathcal{V}} dim(V)$,
then
\[|\mathcal{V}| \leq \qbinom{n}{s}.\]
Moreover, if $|\mathcal{V}| = \qbinom{n}{s}$, then $\mathcal{L} = \{0, 1, \dots, s - 1\}$.
\end{theorem}

We remark that the bound in Theorem \ref{thm1.16} is best possible by taking the example of the family of all $s$-dimensional subspaces of an $n$-dimensional vector space with $\mathcal{L}=\{0, 1, \dots, s - 1\}$.

Note that if $\mathcal{V}$ is an $\mathcal{L}$-intersecting (i.e. $dim(V_{i} \cap V_{j}) \in \mathcal{L}$ for any distinct subspaces $V_{i}$ and $V_{j}$ in $\mathcal{V}$) family of subspaces of an $n$-dimensional vector space satisfying $dim(V_{i}) \not \in \mathcal{L}$ for each $V_{i} \in \mathcal{V}$, then $\mathcal{V}$ is Sperner. Theorem \ref{thm1.16} implies immediately the next consequence which can be viewed as a vector space analogue of Theorem \ref{thm1.9}.

\begin{corollary}\label{cor1.17}
Let $\mathcal{L}=\{l_1, l_2, \dots, l_s\}$ be a set of $s$ nonnegative integers with $0 \leq l_{1} < l_{2} < \cdots < l_{s}$ and $\mathbb{F}_{q}$ be a finite field of order $q$. Suppose that $\mathcal{V}$ is a family of subspaces of an $n$-dimensional vector space over $\mathbb{F}_{q}$ satisfying that (i) $dim(V_{i} \cap V_{j}) \in \mathcal{L}$ for any distinct subspaces $V_{i}$ and $V_{j}$ in $\mathcal{V}$ and (ii) $dim(V_{i}) \not \in \mathcal{L}$ for each $V_{i} \in \mathcal{V}$.
If $n \geq \max\bigg\{\log_{q}\left((q^{s}-1)
\left[ k^{2} \atop l_{1} + 1   \right] + 1 \right)  + l_{1}, \mbox{ } 2s + 1\bigg\}$, where $k = \max_{V \in \mathcal{V}} dim(V)$,
then
\[|\mathcal{V}| \leq \qbinom{n}{s}.\]
Moreover, if $|\mathcal{V}| = \qbinom{n}{s}$, then $\mathcal{L} = \{0, 1, \dots, s - 1\}$.
\end{corollary}

\section{Proofs of Theorems 1.9 and 1.10}
Throughout this paper, we use ${{[n]} \choose {k}}$ to denote the set of all $k$-subsets of $[n] = \{1, 2, \dots, n\}$.
A family $\mathcal{F}$ of sets is said to be an $H_{k}$-\textit{family} ($k \geq 1$) if for every $\mathcal{G} \subseteq \mathcal{F}$,
$\cap_{G \in \mathcal{G}}G = \emptyset$ implies that $\cap_{G \in \mathcal{G}'}G = \emptyset$ for some $\mathcal{G}'  \subseteq \mathcal{G}$
with $|\mathcal{G}'| \leq k$.
The following lemma is Theorem 1(i) from \cite{bd}.

\begin{lemma}\label{lem2.1}
If $\mathcal{F} \subseteq \cup_{i = 0}^{k}{{[n]} \choose {i}}$, then $\mathcal{F}$ is an $H_{d}$-family for every $d \geq k+1$,
that is, if~ $\cap_{F \in \mathcal{F}}F = \emptyset$, then $\cap_{F \in \mathcal{F}'}F = \emptyset$ for some $\mathcal{F}' \subseteq \mathcal{F}$
with $|\mathcal{F}'| \leq k+1$.
\end{lemma}

The next two lemmas are Lemmas 2.2 and 2.3 from \cite{h1}.

\begin{lemma}\label{lem2.2}
Let $l_{1}$ be a positive integer. Let $\mathcal{H}$ be a family of subsets of $[n]$. Suppose that
$\cap_{H \in \mathcal{H}}H = \emptyset$.
Let $F \subseteq [n]$ be a subset such that $F \not \in \mathcal{H}$ and $|F \cap H| \geq l_{1}$ for each $H \in \mathcal{H}$. Set
$Q = \cup_{H \in \mathcal{H}}H$. Then
\[|Q \cap F| \geq l_{1} + 1.\]
\end{lemma}

\begin{lemma}\label{lem2.3}
Let $\mathcal{H}$ be a family of subsets of $[n]$. Suppose that $t = |\mathcal{H}| \geq 2$ and $\mathcal{H}$ is a $k$-uniform intersecting family. Then
\[|\cup_{H \in \mathcal{H}}H| \leq k + (t - 1)(k - 1).\]
\end{lemma}

Lemma \ref{lem2.3} has the following easy corollary by extending every subset to a $k$-subset in an arbitrary way.

\begin{corollary}\label{cor2.4}
Let $\mathcal{H}$ be a family of subsets of $[n]$ with maximum subset size $k$. Suppose that $t = |\mathcal{H}| \geq 2$ and $\mathcal{H}$ is an intersecting family. Then
\[|\cup_{H \in \mathcal{H}}H| \leq k + (t - 1)(k - 1).\]
\end{corollary}

The next fact is easy to verify.

\begin{lemma}\label{lem2.5}
If $n \geq s^{2}$, then
\[{{n - 2} \choose {s}} + {{n - 2} \choose {s - 1}} + \cdots + {{n - 2} \choose {0}} \leq {{n} \choose {s}}.\]
\end{lemma}

\noindent {\bf Proof.}
Recall that
\[{{n} \choose {s}} = {{n-1} \choose {s}} + {{n-1} \choose {s-1}} = {{n-2} \choose {s}} + {{n-2} \choose {s-1}} + {{n-1} \choose {s-1}}.\]
Since $n \geq s^{2}$, we have
\[{{n-1} \choose {s-1}} = \frac{n-1}{s-1} {{n-2} \choose {s-2}} \geq (s - 1){{n-2} \choose {s-2}} \geq {{n - 2} \choose {s - 2}} + \cdots + {{n - 2} \choose {0}}.\]
Thus, it follows that
\[ \hspace{36mm} {{n - 2} \choose {s}} + {{n - 2} \choose {s - 1}} + \cdots + {{n - 2} \choose {0}} \leq {{n} \choose {s}}. \hspace{36mm} \hfill \Box\]

\vspace{3mm}

We also need the following theorem which is Corollary 1.8 in \cite{lzx}.

\begin{theorem}\label{thm2.6} (Liu, Zhang, and Xiao \cite{lzx}).
Let $\mathcal{L}=\{l_1, l_2, \dots, l_s\}$ be a set of $s$ nonnegative integers with $l_{1} < l_{2} < \cdots < l_{s}$. Suppose that $\mathcal{A} = \{A_1, A_2, \dots, A_{m}\}$ is an $\mathcal{L}$-intersecting family of subsets of $[n]$.
If $n \geq {{k^{2}} \choose {l_{1}+1}}s + l_{1}$, then
\[|\mathcal{A}|\leq {{n - l_{1}} \choose {s}} + {{n - l_{1}} \choose {s - 1}} + \cdots + {{n - l_{1}} \choose {0}}.\]
\end{theorem}

To prove the next lemma,
we will use ${\bf x} = (x_1, x_2, \dots, x_n)$ to denote a vector of $n$ variables with each variable $x_j$ taking values $0$ or $1$. A polynomial $f({\bf x})$ in $n$ variables $x_{i}$, $1 \leq i \leq n$, is called $multilinear$ if the power of each variable $x_{i}$ in each term is at most one. Clearly, if each variable
$x_i$ takes only the values $0$ or $1$, then any polynomial in $n$ variables $x_{i}$, $1 \leq i \leq n$, is multilinear. For a subset $F$ of $[n]$, we define the characteristic vector ${\bf v}_{F}$ of $F$ to be the vector ${\bf v}_{F} = (v_1, v_2, \dots, v_{n})$ with $v_{j} = 1 $ if $j \in F$ and $v_{j} = 0$ otherwise.

\begin{lemma}\label{lem2.7}
Let $\mathcal{L}=\{l_1, l_2, \dots, l_s\}$ be a set of $s$ nonnegative integers with $l_{1} < l_{2} < \cdots < l_{s}$. Suppose that $\mathcal{A} = \{A_1, A_2, \dots, A_{m}\}$ and $\mathcal{B} = \{B_1, B_2, \dots, B_{m}\}$ are two families of subsets of $[n]$ such that (i) $|A_{i} \cap B_{j}| \in \mathcal{L}$ for any $i \neq j$ and
(ii) $A_{i} \subseteq B_{i}$  and $|A_{i}| \not \in \mathcal{L}$ for every $1 \leq i \leq m$.
Then
\[m = |\mathcal{A}|\leq {{n - 1} \choose {s}} + {{n - 1} \choose {s - 1}} + \cdots + {{n - 1} \choose {0}}.\]
\end{lemma}

\noindent {\bf Proof.}
Suppose that $\mathcal{A} = \{A_1, A_2, \dots, A_{m}\}$ and
$\mathcal{B} = \{B_{1}, B_{2}, \dots, B_{m}\}$ are two families of subsets of $[n]$  satisfying the conditions in the lemma.
Without loss of generality, we may assume that $n \in A_{j}$ for $1 \leq j \leq t$ and $n \not \in A_{j} $  for $ j \geq t+1$. Then $n \in B_{j} $  for $1 \leq j \leq t$.

For each $A_{j} \in \mathcal{A}$, define
\[f_{A_{j}}({\bf x})=\prod_{i = 1}^{s}({\bf v}_{A_j}\cdot {\bf x} - l_{i}),\]
where ${\bf x} = (x_1, x_2, \dots, x_n)$ is a vector of $n$ variables with each variable $x_i$ taking values $0$ or $1$, and ${\bf x} \cdot {\bf y}$ is the inner product of vectors ${\bf x}$ and ${\bf y}$.
Then each $f_{A_{j}}({\bf x})$ is a multilinear polynomial of degree at most $s$. Since ${\bf v}_{A_i} \cdot {\bf v}_{B_j} = |A_{i} \cap B_{j}|$ for any pair $i, j$,  condition (i) implies $f_{A_{j}}({\bf v}_{B_i}) = 0$ for $i \neq j$; and
condition (ii) implies $f_{A_{j}}({\bf v}_{B_j}) \neq 0$ for each $1 \leq j \leq m$.

Let $Q$ be the family of subsets of $[n]$ containing $n$ with sizes at most $s$. Then $|Q| = \sum_{i = 0}^{s-1}{{n-1} \choose {i}}$.
For each $R \in Q$, define
\[g_{R}({\bf x}) = (1-x_{n})\prod_{j \in R, j \neq n}x_{j}.\]
Then each $g_{R}({\bf x})$ is a multilinear polynomial of degree at most $s$.

We now proceed to show that the polynomials in
\[\{f_{A_i}({\bf x}) \mid 1 \leq i \leq m\} \cup \{g_{R}({\bf x}) \mid R \in Q\} \]
are linearly independent.
Suppose that we have a linear combination of these polynomials that equals to zero:
\begin{equation}
 \sum_{i = 1}^{m} \alpha_{i}f_{A_i}({\bf x})+ \sum_{R \in Q}\beta_{R}g_{R}({\bf x}) = 0.
\label{2.2}\end{equation}
We show that all coefficients must be zero as follows.

\noindent {\bf Claim 1.} $\alpha_{j} = 0$ for each $1 \leq j \leq t$ ($n \in A_{j} \subseteq B_{j}$ for each $j \leq t$).

Suppose, to the contrary, that $c \leq t$ is the smallest integer such that $\alpha_{c} \neq 0$. Since $n \in A_{c} \subseteq B_{c}$, $g_{R}({\bf v}_{B_c}) = 0$ for every $R \in Q$.
Recall that $f_{A_{j}}({\bf v}_{B_c}) = 0$ for each $j \neq c$ and $f_{A_{c}}({\bf v}_{B_c}) \neq 0$.
By evaluating (2) with ${\bf x} = {\bf v}_{B_c}$, we obtain
$\alpha_{c}f_{A_c}({\bf v}_{B_c}) = 0$ which implies that $\alpha_{c} = 0$, a contradiction. Thus, Claim 1 holds.

\noindent {\bf Claim 2.} $\alpha_{j} = 0$ for each $j \geq t + 1$ ($n \not \in A_{j}$ for $j \geq t + 1$).

By Claim 1 and Equation (2), we have
\begin{equation}
\sum_{i = t + 1}^{m} \alpha_{i}f_{A_i}({\bf x})+ \sum_{R \in Q}\beta_{R}g_{R}({\bf x}) = 0,
\label{2.3}\end{equation}
Suppose, to the contrary, that $c \geq t + 1$ is the smallest integer such that $\alpha_{c} \neq 0$. Since $n \not \in A_{c}$, the $n$-th coordinate of
${\bf v}_{A_c}$ is $0$. For each $B_{j}$ with $j \geq t + 1$,
let ${\bf v}'_{B_{j}} = {\bf v}_{B_{j}} + (0, \dots, 0, 1)$ if $n \not \in B_{j}$ and ${\bf v}'_{B_{j}} = {\bf v}_{B_{j}}$ if $n \in B_{j}$. Then
$g_{R}({\bf v}'_{B_c}) = 0$ for every $R \in Q$.  For each $j \geq t + 1$, since $n \not \in A_{j} $, $f_{A_j}({\bf v}'_{B_c}) = f_{A_j}({\bf v}_{B_c})$. Thus, by evaluating (3) with ${\bf x} = {\bf v}'_{B_c}$, we have
$\alpha_{c}f_{A_c}({\bf v}_{B_c}) = \alpha_{c}f_{A_c}({\bf v}'_{B_c}) = 0$ which implies that $\alpha_{c} = 0$, a contradiction. Thus, Claim 2 is true.

By Claim 2 and Equation (3), we now have
\begin{equation}
\sum_{R \in Q}\beta_{R}g_{R}({\bf x}) = 0.
\label{2.4}\end{equation}
Setting $x_{n} = 0$ in Equation (4), we have $\sum_{R \in Q}\beta_{R}g'_{R}({\bf x}) = 0$, where
$g'_{R}({\bf x}) = \prod_{j \in R, j \neq n}x_{j}$. It is easy to see that the polynomials $g'_{R}({\bf x})$, $R \in Q$, are linearly independent. Thus, we have $\beta_{R} = 0$ for all $R \in Q$.

We have shown that the polynomials in
\[\{f_{A_i}({\bf x}) \mid 1 \leq i \leq m\} \cup \{g_{R}({\bf x}) \mid R \in Q\} \]
are linearly independent. Since the dimension of the vector space of multilinear polynomials of degree at most $s$ is
\[ {{n}\choose {s}}+ {{n} \choose {s-1}}+
\cdots + {{n} \choose {0}},\]
it follows that
\[m + |Q| \leq {{n}\choose {s}}+ {{n} \choose {s-1}}+
\cdots + {{n} \choose {0}},\]
which implies that
\[  m = |\mathcal{A}| \leq {{n-1}\choose {s}}+ {{n-1} \choose {s-1}}+ \cdots + {{n-1} \choose {0}}.  \]
$\hfill \Box$

\vspace{3mm}

As a consequence of Lemma \ref{lem2.7} by taking $\mathcal{A} = \{A_1, A_2, \dots, A_{m}\} = \mathcal{B} = \{B_1, B_2, \dots, B_{m}\}$, one has the next immediate corollary.

\begin{corollary}\label{cor2.8}
Let $\mathcal{L}=\{l_1, l_2, \dots, l_s\}$ be a set of $s$ nonnegative integers with $0 \leq l_{1} < l_{2} < \cdots < l_{s}$. Suppose that $\mathcal{A} = \{A_1, A_2, \dots, A_{m}\}$ is a family of subsets of $[n]$ such that (i) $|A_{i} \cap A_{j}| \in \mathcal{L}$ for any $i \neq j$ and
(ii) $|A_{i}| \not \in \mathcal{L}$ for every $1 \leq i \leq m$.
Then
\[|\mathcal{A}|\leq {{n - 1} \choose {s}} + {{n - 1} \choose {s - 1}} + \cdots + {{n - 1} \choose {0}}.\]
\end{corollary}

We are now ready to provide a proof for Theorem 1.9.

\vspace{3mm}

\noindent {\bf Proof of Theorem 1.9.}
Let $\mathcal{L} = \{l_1, l_2, \dots, l_s\}$ with $0 \leq l_{1} < l_{2} < \cdots < l_{s}$.
We prove by induction on $s \geq 1$. For $s = 1$, it follows from Corollary \ref{cor2.8} that
\[ m = |\mathcal{A}| \leq {{n-1}\choose {1}}+ {{n-1} \choose {0}} = {{n}\choose {1}}.  \]
Assume that the result holds for any positive integer $s = r$. We now consider that $s = r + 1 \geq 2$.
We consider the following two cases.

\noindent {\bf Case 1.} $0 \in \mathcal{L}$, i.e., $l_{1} = 0$. If $\mathcal{L}=\{0, 1, \dots, s - 1\}$, then the result follows from Theorem \ref{thm1.8}.
Assume that $\mathcal{L} \neq \{0, 1, \dots, s - 1\}$.
Let $\mathcal{A}' = \{A_{j} \in \mathcal{A} \mid  |A_{j}| \leq s\}$ and $\mathcal{L}' = \mathcal{L} \cap \{0, 1, 2, \dots, s - 1\}$. Then $|\mathcal{L}'| \leq s - 1$.
Clearly, $\mathcal{A}'$ is an $\mathcal{L}'$-intersecting family satisfying $|A_{j}| \not \in \mathcal{L}'$ for every $A_{j} \in \mathcal{A}'$. By the induction hypothesis, we have
\[|\mathcal{A}'| \leq {{n} \choose {s - 1}}.\]
For each $h \in [n]$, let
\[\mathcal{A}_{h} = \{A_{j} \in \mathcal{A} \mid  h \in A_{j}\},\]
\[\mathcal{A}'_{h} = \{A_{j}\setminus \{h\} \mid  A_{j} \in \mathcal{A}_{h}\}.\]
Then $\mathcal{A}'_{h}$ is an $\mathcal{L}^{*}$-intersecting family on the set $[n] \setminus \{h\}$ satisfying $|A_{j}\setminus \{h\}| \not \in \mathcal{L}^{*}$ for each $A_{j} \in \mathcal{A}_{h}$,
where $\mathcal{L}^{*} = \{l_{2} - 1, \dots, l_{s} - 1\}$.
By Corollary \ref{cor2.8}, we have for each $h \in [n]$,
\[|\mathcal{A}_{h}| = |\mathcal{A}'_{h}| \leq {{n - 1} \choose {s - 1}} + {{n - 1} \choose {s - 2}} + \cdots + {{n - 1} \choose {0}}.\]
Note that for each $A_{j} \in \mathcal{A} \setminus \mathcal{A}'$, $|A_{j}| \geq s + 1$. It follows that
\[(m - |\mathcal{A}'|)(s + 1) \leq \sum_{h \in [n]}|\mathcal{A}_{h}| \leq n \bigg[{{n - 1} \choose {s - 1}} + {{n - 1} \choose {s - 2}} + \cdots + {{n - 1} \choose {0}}\bigg].\]
Note that $k = \max_{1 \leq j \leq m} |A_{j}| \geq s + 1$ as $\mathcal{L} = \{l_1, l_2, \dots, l_s\} \neq \{0, 1, \dots, s - 1\}$. We have $n \geq {{k^{2}} \choose {l_{1}+1}}s + l_{1} \geq s^{3} + 1$. It follows that
\[m = |\mathcal{A}| = |\mathcal{A}'| + (m - |\mathcal{A}'|) \leq {{n} \choose {s - 1}} + \frac{n}{s + 1}\sum_{j = 0}^{s - 1}{{n - 1} \choose {j}}
\leq  {{n} \choose {s}},\]
where the last inequality holds because
\[{{n} \choose {s}} - {{n} \choose {s - 1}} = {{n - 1} \choose {s}} - {{n - 2} \choose {s - 2}} \geq
\frac{n}{s + 1} \bigg[{{n - 1} \choose {s -1}} + (s -1){{n - 1} \choose {s - 2}}\bigg] \geq \frac{n}{s + 1}\sum_{j =0}^{s-1}{{n - 1} \choose {j}}.\]

\noindent {\bf Case 2.} $0 \not \in \mathcal{L}$, i.e., $l_{1} > 0$.
We consider two subcases.

\noindent {\bf Subcase 2.1.} $\cap_{A_{i} \in \mathcal{A}}A_{i} = \emptyset$.

By Lemma \ref{lem2.1}, there exists a subfamily $\mathcal{A}' \subseteq \mathcal{A}$ with $|\mathcal{A}'| = k + 1$ such that
$\cap_{A_{i} \in \mathcal{A}'}A_{i} = \emptyset$.
Let $M = \cup_{A_{i} \in \mathcal{A}'}A_{i}$. Then $|M| \leq k + k(k - 1) = k^{2}$ by Corollary \ref{cor2.4}.
Since $\mathcal{A}$ is an $\mathcal{L}$-intersecting family satisfying $|A_{i}| \not \in \mathcal{L}$ for every $1 \leq i \leq m$, $|A_{i}| \geq l_{1} + 1$ for every $A_{i} \in \mathcal{A}$. It follows from Lemma \ref{lem2.2} that
\begin{equation}
|M \cap A_{i}| \geq l_{1} + 1 \mbox{ for all } A_{i} \in \mathcal{A}.
\label{2.1}\end{equation}
Let $T$ be a given subset of $M$ such that $|T| = l_{1} + 1$. Define
\[\mathcal{A}(T) = \{A_{i} \in \mathcal{A} \mid T \subseteq M \cap A_{i}\}.\]
Set $\mathcal{L}' = \{l_{2}, l_{3}, \dots, l_{s}\}$. Then $|\mathcal{L}'| = s - 1$.
Let
\[\mathcal{A}'(T) = \{A_{i}\setminus T \mid A_{i} \in  \mathcal{A}(T)\}.\]
Since $\mathcal{A}$ is $\mathcal{L}$-intersecting and $|A_{i} \cap A_{j}| \geq |T| = l_{1} + 1$,
$\mathcal{A}(T)$ is $\mathcal{L}'$-intersecting. Moreover, $\mathcal{A}'(T)$ is an $\mathcal{L}^{*}$-intersecting family such that $|A_{i}\setminus T| \not \in \mathcal{L}^{*}$ for every $A_{i} \in  \mathcal{A}(T)$, where
$\mathcal{L}^{*} = \{l_{2} - l_{1} - 1, l_{3} - l_{1} - 1, \dots, l_{s} - l_{1} - 1\}$.
By the induction hypothesis,
we have
\[|\mathcal{A}'(T)| \leq {{n - l_{1} - 1} \choose {s - 1}}.\]
By (5), we have
\[\mathcal{A} = \cup_{T \subseteq M, \mbox{ } |T| = l_{1} + 1}\mathcal{A}(T).\]
It follows that
\[|\mathcal{A}| \leq {{k^{2}} \choose {l_{1}+1}}{{n - l_{1} - 1} \choose {s - 1}} = {{k^{2}} \choose {l_{1}+1}}\frac{s}{n - l_{1}}{{n - l_{1}} \choose {s}}
\leq {{n - l_{1}} \choose {s}} \leq {{n} \choose {s}}.\]

\noindent {\bf Subcase 2.2.} $\cap_{A_{i} \in \mathcal{A}}A_{i} \neq \emptyset$.
Let $Z = \cap_{A_{i} \in \mathcal{A}}A_{i}$ and $t = |Z|$. We may assume that $1 \leq t \leq l_{1}$.
If $l_{1} \geq 2$, then it follows from Theorem \ref{thm2.6} and Lemma \ref{lem2.5} that
$$
\begin{aligned}
 |\mathcal{A}|&\leq {{n - l_{1}} \choose {s}} + {{n - l_{1}} \choose {s - 1}} + \cdots + {{n - l_{1}} \choose {0}}\\
&\leq {{n - 2} \choose {s}} + {{n - 2} \choose {s - 1}} + \cdots + {{n - 2} \choose {0}}\\
& \leq {{n} \choose {s}}.
\end{aligned}
$$
Assume that $l_{1} = 1$. Then $|Z| = t = l_{1} = 1$. Let $\mathcal{L}'' = \{0, l_{2} - 1, \dots, l_{s} - 1\}$ and
\[\mathcal{A}[Z] = \{A_{i} \setminus Z \mid  A_{i} \in \mathcal{A}\}.\]
Then it is clear that $\mathcal{A}[Z]$ is an $\mathcal{L}''$-intersecting family of subsets on $[n] \setminus Z$ satisfying
$|A_{i} \setminus Z| \not \in \mathcal{L}''$ for every $A_{i} \setminus Z  \in \mathcal{A}[Z]$.
Since $0 \in \mathcal{L}''$, it follows from Case 1 that
\[|\mathcal{A}| = |\mathcal{A}[Z]| \leq {{n} \choose {s}}.\]
Therefore, the theorem follows.
\hfill$\Box$

\begin{lemma}\label{lem2.9}
Let $\mathcal{L}=\{l_1, l_2, \dots, l_s\}$ be a set of $s$ nonnegative integers with $0 \leq l_{1} < l_{2} < \cdots < l_{s}$. Suppose that $\mathcal{A} = \{A_1, A_2, \dots, A_{m}\}$ and $\mathcal{B} = \{B_1, B_2, \dots, B_{m}\}$ are families of subsets of $[n]$ such that (i) $|A_{i} \cap B_{j}| \in \mathcal{L}$ for any $i \neq j$;
(ii) $A_{i} \subseteq B_{i}$  and $|A_{i}| \not \in \mathcal{L}$ for every $1 \leq i \leq m$; and (iii)
$A_{i} = B_{i}$ for $i \leq k + 1$ and $\cap_{1 \leq j \leq k + 1}A_{j} = \cap_{A_{j} \in \mathcal{A}}A_{j}$, where $k$ is the maximum subset size in $\mathcal{B}$.
If  $n \geq {{k^{2}} \choose {l_{1}+1}}s + l_{1}$, then
\[|\mathcal{A}|\leq {{n} \choose {s}}.\]
\end{lemma}

\noindent {\bf Proof.}
We prove by induction on $s$. For $s =1$, it follows from Lemma \ref{lem2.7} that
\[ m = |\mathcal{A}| \leq {{n-1}\choose {1}}+ {{n-1} \choose {0}} = {{n}\choose {1}}.  \]
Assume that the result holds for any positive integer $s = r$. We now consider that $s = r + 1 \geq 2$.
We consider the following two cases.

\noindent {\bf Case 1.} $\cap_{A_{j} \in \mathcal{A}}A_{j} \neq \emptyset$. Let $z \in \cap_{A_{j} \in \mathcal{A}}A_{j} \subseteq \cap_{A_{j} \in \mathcal{A}}B_{j}$.
Set
\[\mathcal{A}(z) = \{A_{j}\setminus \{z\} \mid  A_{j} \in \mathcal{A}\},\]
\[\mathcal{B}(z) = \{B_{j}\setminus \{z\} \mid  B_{j} \in \mathcal{B}\}.\]
Then $\mathcal{A}(z)$ and $\mathcal{B}(z)$ are cross $\mathcal{L}'$-intersecting and $|A_{i} \setminus \{z\}| \not \in \mathcal{L}'$ for each $i \leq m$,
where $\mathcal{L}' = \{l_{1} - 1, \dots, l_{s} - 1\}$. It follows from Lemmas \ref{lem2.5} and \ref{lem2.7} that
\[|\mathcal{A}|\leq {{n - 1 - 1} \choose {s}} + {{n - 1 - 1} \choose {s - 1}} + \cdots + {{n - 1 - 1} \choose {0}} \leq {{n} \choose {s}}.\]

\noindent {\bf Case 2.} $\cap_{A_{j} \in \mathcal{A}}A_{j} = \emptyset$.
Then $\cap_{1 \leq j \leq k + 1}A_{j} = \cap_{1 \leq j \leq k + 1}B_{j} = \emptyset$ by assumption (iii).
We consider two subcases.

\noindent {\bf Subcase 2.1.} $l_{1} = 0$. i.e., $0 \in \mathcal{L}$.
Suppose that $\mathcal{L} = \{0, 1, 2, \dots, s - 1\}$. In this case, by assumptions (i) and (ii), it follows that  $|A_{j}| \geq s$ for each $1 \leq j \leq m$, and for every subset $S \subseteq [n]$ with $|S| = s$, there exists at most one subset $A_{j} \in \mathcal{A}$ such that $S \subseteq A_{j}$. Define $\mathcal{A}' = \{A'_{j} \mid  A'_{j} \subseteq A_{j} \in \mathcal{A} \mbox{ with } |A'_{j}| = s\}$. Then $\mathcal{A}'$ is $\mathcal{L}$-intersecting since
$A'_{i} \cap A'_{j} \subseteq A_{i} \cap A_{j} \subseteq A_{i} \cap B_{j}$ and $|A_{i} \cap B_{j}| \in \mathcal{L} = \{0, 1, 2, \dots, s - 1\}$ for any $i \neq j$.
Thus, $\mathcal{A}'$ is $s$-uniform $\mathcal{L}$-intersecting and it follows from Theorem \ref{thm1.2} that
\[m = |\mathcal{A}| = |\mathcal{A}'| \leq {{n} \choose {s}}.\]
Now assume that $\mathcal{L} \neq \{0, 1, 2, \dots, s - 1\}$.
Let $\mathcal{A}' = \{A_{j} \in \mathcal{A} \mid |A_{j}| \leq s\}$ and $\mathcal{L}' = \mathcal{L} \cap \{0, 1, 2, \dots, s - 1\}$. Then $|\mathcal{L}'| \leq s - 1$.
Let $\mathcal{B}' = \{B_{j} \in \mathcal{B} \mid A_{j} \in \mathcal{A}'\}$. Clearly, $\mathcal{A}'$ and $\mathcal{B}'$ are cross $\mathcal{L}'$-intersecting. By the induction hypothesis, we have
\[|\mathcal{A}'| \leq {{n} \choose {s - 1}}.\]
For each $h \in [n]$, let
\[\mathcal{A}_{h} = \{A_{j} \in \mathcal{A} \mid  h \in A_{j}\}\]
and let
\[\mathcal{B}_{h} = \{B_{j} \in \mathcal{B} \mid A_{j} \in \mathcal{A}_{h}\}.\]
Then $\mathcal{A}_{h}$ and $\mathcal{B}_{h}$ are cross $\mathcal{L}^{*}$-intersecting, where $\mathcal{L}^{*} = \mathcal{L}\setminus \{0\}$.
By Lemma \ref{lem2.7}, we have, for each $h \in [n]$,
\[|\mathcal{A}_{h}| \leq {{n - 1} \choose {s - 1}} + {{n - 1} \choose {s - 2}} + \cdots + {{n - 1} \choose {0}}.\]
Note that for each $A_{j} \in \mathcal{A} \setminus \mathcal{A}'$, $|A_{j}| \geq s + 1$. It follows that
\[(m - |\mathcal{A}'|)(s + 1) \leq \sum_{h \in [n]}|\mathcal{A}_{h}| \leq n \left[{{n - 1} \choose {s - 1}} + {{n - 1} \choose {s - 2}} + \cdots + {{n - 1} \choose {0}}\right].\]
Since $n \geq {{k^{2}} \choose {l_{1}+1}}s + l_{1} \geq s^{3} + 1$, we have
\[m = |\mathcal{A}| = |\mathcal{A}'| + (m - |\mathcal{A}'|) \leq {{n} \choose {s - 1}} + \frac{n}{s + 1}\sum_{j = 0}^{s - 1}{{n - 1} \choose {j}}
\leq  {{n} \choose {s}}.\]

\noindent {\bf Subcase 2.2.} $l_{1} > 0$.
Let $M = \cup_{1 \leq i \leq k + 1}B_{i}$. Then $|M| \leq k + k(k - 1) = k^{2}$ by Corollary \ref{cor2.4}.
Since $\mathcal{A}$ and $\mathcal{B}$ are cross $\mathcal{L}$-intersecting and $|A_{i}| \not \in \mathcal{L}$ for every $1 \leq i \leq m$, $|A_{i}| \geq l_{1} + 1$ for every $1 \leq i \leq m$. Note that $A_{i} \subseteq B_{i} \subseteq M$ for each $1 \leq i \leq k + 1$. It follows from Lemma \ref{lem2.2} that
\begin{equation}
|M \cap A_{i}| \geq l_{1} + 1 \mbox{ for all } A_{i} \in \mathcal{A}.
\label{2.5}\end{equation}
Let $T$ be a given subset of $M$ such that $|T| = l_{1} + 1$. Define
\[\mathcal{A}(T) = \{A_{i} \in \mathcal{A} \mid T \subseteq M \cap A_{i}\}.\]
Then
\[\mathcal{B}(T) = \{B_{i} \in \mathcal{B} \mid T \subseteq M \cap B_{i}\} = \{B_{i} \in \mathcal{B} \mid A_{i} \in \mathcal{A}(T)\}.\]
Let
\[\mathcal{A}'(T) = \{A_{i}\setminus T \mid A_{i} \in  \mathcal{A}(T)\},\]
\[\mathcal{B}'(T) = \{B_{i}\setminus T \mid B_{i} \in  \mathcal{B}(T)\}.\]
Since $\mathcal{A}$ and $\mathcal{B}$ are cross $\mathcal{L}$-intersecting and $|A_{i} \cap B_{j}| \geq |T| = l_{1} + 1$ for every $A_{i} \in \mathcal{A}(T)$ and $B_{j} \in \mathcal{B}(T)$,
$\mathcal{A}(T)$ and $\mathcal{B}(T)$ are cross $\mathcal{L}'$-intersecting, where $\mathcal{L}' = \{l_{2}, l_{3}, \dots, l_{s}\}$ and $|\mathcal{L}'| = s - 1$.
Thus, $\mathcal{A}'(T)$ and $\mathcal{B}'(T)$ are cross $\mathcal{L}^{*}$-intersecting families on $[n - l_{1} - 1]$, where
$\mathcal{L}^{*} = \{l_{2} - l_{1} - 1, l_{3} - l_{1} - 1, \dots, l_{s} - l_{1} - 1\}$.
By the induction hypothesis,
we have
\[|\mathcal{A}'(T)| \leq {{n - l_{1} - 1} \choose {s - 1}}.\]
By (6), we have
\[\mathcal{A} = \cup_{T \subseteq M, \mbox{ } |T| = l_{1} + 1}\mathcal{A}(T).\]
It follows that
\[\hspace{12mm} |\mathcal{A}| \leq {{k^{2}} \choose {l_{1}+1}}{{n - l_{1} - 1} \choose {s - 1}}
= {{k^{2}} \choose {l_{1}+1}}\frac{s}{n - l_{1}}{{n - l_{1}} \choose {s}}
\leq {{n - l_{1}} \choose {s}} \leq {{n} \choose {s}}. \hspace{12mm} \hfill \Box\]

\vspace{3mm}

We now provide our proof for Theorem \ref{thm1.10}.

\vspace{3mm}

\noindent {\bf Proof of Theorem 1.10.}
We prove by induction on $t \geq 2$. For $t = 2$, the result follows from Theorem \ref{thm1.9}.
Assume that the result holds for $(t - 1)$-wise $\mathcal{L}$-intersecting families. We consider a $t$-wise $\mathcal{L}$-intersecting family
$\mathcal{A} = \{A_1, A_2, \dots, A_{m}\}$.
We partition $\mathcal{A}$ into two families $\mathcal{B}$ and $\mathcal{F}$ satisfying the following properties: $\mathcal{F}$ is $(t - 1)$-wise $\mathcal{L}$-intersecting and there exists a family $\mathcal{C}$ such that the pair $(\mathcal{B}, \mathcal{C})$ satisfies the assumptions in Lemma \ref{lem2.9}.

To obtain the desired partition, we first construct families $\mathcal{B}$, $\mathcal{C}$ and $\mathcal{F}$ by repeating the following procedure: By Lemma \ref{lem2.1},
without loss of generality, assume $\cap_{1 \leq j \leq k + 1}A_{j} = \cap_{A_{j} \in \mathcal{A}}A_{j}$.
Set $B_{i} = C_{i} = A_{i}$ for $i \leq k + 1$.
For every $k + 1 \leq d \leq m - 1$, suppose that after step $d$, we have constructed families $\mathcal{B} = \{B_{1}, B_{2}, \dots, B_{d}\}$ and
$\mathcal{C} = \{C_{1}, C_{2}, \dots, C_{d}\}$. At step $d + 1$, suppose that there are indices $i_{1} < i_{2} < \cdots < i_{t - 1}$ such that
$|A_{i_{1}} \cap A_{i_{2}} \cap \cdots \cap A_{i_{t-1}}| \not \in \mathcal{L}$. With relabeling if necessary, we assume $i_{1} = d + 1$.
Set $B_{d+1} = A_{i_{1}} = A_{d+1}$ and $C_{d+1} = A_{i_{1}} \cap A_{i_{2}} \cap \cdots \cap A_{i_{t-1}}$. Then $C_{d+1} \subseteq B_{d+1}$,
$|B_{d+1} \cap C_{d+1}| = |C_{d+1}| \not \in \mathcal{L}$, and $|B_{i} \cap C_{j}| \in \mathcal{L}$ for $1 \leq i \neq j \leq d+1$. Update $d$ by $d+1$ and proceed to the next step. Continue this process until we can not proceed further. Then set $\mathcal{F} = \mathcal{A} \setminus \mathcal{B}$.

Then, by the construction, we have $m = |\mathcal{B}| + |\mathcal{F}|$, the families $\mathcal{B}$ and $\mathcal{C}$ satisfy the assumptions (i)--(iii) in Lemma \ref{lem2.9}, and $\mathcal{F}$ is $(t-1)$-wise $\mathcal{L}$-intersecting. It follows from Lemma \ref{lem2.9} that
\[|\mathcal{B}| \leq {{n} \choose {s}}.\]
By the induction hypothesis, we have
\[ |\mathcal{F}| \leq (t - 2) {{n} \choose {s}}.\]
It follows that
\[m = |\mathcal{B}| + |\mathcal{F}| \leq {{n} \choose {s}} + (t - 2) {{n} \choose {s}} = (t - 1) {{n} \choose {s}}.\]
$\hfill \Box$

\section{Proofs of Theorems 1.15 and 1.16}
Note that for $\qbinom{n}{k}$ defined by (1), one has
\[\qbinom{n}{k} = q^{k}\qbinom{n-1}{k} + \qbinom{n-1}{k-1}.\]
First, we show that $\left[n \atop k \right]$ has the following property similar to the binomial coefficients.

\begin{lemma}\label{lem3.1}
For $0 \leq k < l \leq \lfloor \frac{n}{2} \rfloor$, we have
\[{\rm (i)} \mbox{ } \qbinom{n}{k} = \qbinom{n}{n-k}, \hspace{120mm}\]
\[{\rm (ii)} \mbox{ } \qbinom{n}{k} < \qbinom{n}{l}.  \hspace{126mm}\]
\end{lemma}

\noindent {\bf Proof.}
For (i), it follows from the fact
\[\qbinom{n}{n-k} = \prod_{0 \leq i \leq n-k-1}\frac{q^{n-i}-1}{q^{n-k-i} - 1}
= \frac{(q^{n}-1)(q^{n-1}-1)\cdots (q^{n-(n-k)+1}-1)}{(q^{n-k}-1)(q^{n-k-1}-1)\cdots (q - 1)} \hspace{25mm}\]
\[ = \frac{(q^{n}-1)(q^{n-1}-1)\cdots (q^{n-(n-k)+1}-1)}{(q^{n-k}-1)(q^{n-k-1}-1)\cdots (q - 1)} \cdot
\frac{(q^{n}-1)(q^{n-1}-1)\cdots (q-1)}{(q^{n}-1)(q^{n-1}-1)\cdots (q-1)}\]
\[ = \frac{(q^{n}-1)(q^{n-1}-1)\cdots (q^{n-k+1}-1)}{(q^{k}-1)(q^{k-1}-1)\cdots (q - 1)} = \qbinom{n}{k}. \hspace{45mm}\]

\noindent For (ii), it suffices to show that for $k + 1 \leq \lfloor \frac{n}{2} \rfloor$,
\[\qbinom{n}{k} < \qbinom{n}{k+1}.\]
Note that
\[\qbinom{n}{k+1} = \frac{q^{n - k} - 1}{q^{k+1} - 1}\qbinom{n}{k}.\]
For $k + 1 \leq \lfloor \frac{n}{2} \rfloor$, $n - k > k + 1$ and it follows that
\[\qbinom{n}{k} < \qbinom{n}{k+1}.\]
$\hfill \Box$

The next fact can be verified easily.

\begin{lemma}\label{lem3.2}
Let $s \geq 1$, $q \geq 2$, and $n \geq 2s + 1$.
Then
\[ \qbinom{n-1}{s}
+  \qbinom{n-1}{s-1} + \cdots +  \qbinom{n-1}{0} <  \qbinom{n}{s}.\]
\end{lemma}

\noindent {\bf Proof.}
Since $q \geq 2$ and $n \geq 2s + 1$, by Lemma \ref{lem3.1}, we have
\[\qbinom{n}{s} = q^{s} \qbinom{n-1}{s} + \qbinom{n-1}{s-1} \hspace{52mm}\]
\[> s\qbinom{n-1}{s} + \qbinom{n-1}{s-1}  \hspace{46mm}\]
\[= \qbinom{n-1}{s} + \qbinom{n-1}{s-1} + (s - 1)\qbinom{n-1}{s}  \hspace{18mm}\]
\[\geq \qbinom{n-1}{s}
+ \qbinom{n-1}{s-1} + \qbinom{n-1}{s-2} + \cdots + \qbinom{n-1}{0}.\hspace{1mm}\]
$\hfill \Box$

\vspace{3mm}

The following LYM inequality for Sperner family of subspaces of an $n$-dimensional vector space $\mathcal{W}$ over a finite field $\mathbb{F}_{q}$
is implied by Theorem 13 in \cite{b}.

\begin{theorem}\label{thm3.3} (Bey, \cite{b}).
Suppose that
$\mathcal{V}$ is a Sperner family of subspaces of an $n$-dimensional vector space $\mathcal{W}$ over a finite field $\mathbb{F}_{q}$.
For $0 \leq k \leq n$, let $\mathcal{V}_{k} = \{V \in \mathcal{V} \mid dim(V) = k\}$.
Then
\[\sum_{k = 0}^{n}\frac{|\mathcal{V}_{k}|}{\qbinom{n}{k}} \leq 1.\]
\end{theorem}

By Lemma \ref{lem3.1} and Theorem \ref{thm3.3}, we have the next consequence immediately.

\begin{theorem}\label{thm3.4}
Let $n \geq 2s + 1$.
Suppose that
$\mathcal{V}$ is a Sperner family of subspaces of an $n$-dimensional vector space $\mathcal{W}$ over a finite field $\mathbb{F}_{q}$ such that
$dim(V) \leq s$ for each $V \in \mathcal{V}$. Then
\[|\mathcal{V}| \leq \qbinom{n}{s} \]
with equality only if $\mathcal{V}$ is the set of all $s$-dimensional subspaces of $\mathcal{W}$.
\end{theorem}

\noindent {\bf Proof of Theorem 1.15.} This theorem follows directly from Theorem 1.15 in \cite{lzx} and Lemma \ref{lem3.2}.
$\hfill \Box$

\vspace{3mm}

\noindent {\bf Proof of Theorem 1.16.}
Suppose that $\mathcal{V} = \{V_{1}, V_{2}, \dots, V_{m}\}$ is a Sperner family of subspaces of an $n$-dimensional vector space $\mathcal{W}$ over $\mathbb{F}_{q}$ satisfying that $dim(V_{i} \cap V_{j}) \in \mathcal{L}$ for any distinct subspaces $V_{i}$ and $V_{j}$ in $\mathcal{V}$.
By Theorem \ref{thm1.15}, we may assume that $0 \in \mathcal{L}$.
We prove by induction on $s \geq 1$. For $s = 1$, we have $\mathcal{L} = \{0\}$ and the result follows from Case 1 below. Assume that the result holds for any positive integer $s = h$. We now consider $s = h + 1 \geq 2$. We consider the following two cases.

\noindent {\bf Case 1.} $\mathcal{L} = \{0, 1, \dots, s-1\}$.

Then for any fixed $h \geq s$, every subspace of dimension $h$ can be a subspace of at most one subspace in $\mathcal{V}$ as $\mathcal{V}$ is $\mathcal{L}$-intersecting.
Construct a new family $\mathcal{V}' = \{V'_{i} \mid 1 \leq i \leq m\}$, where for each $1 \leq i \leq m$, $V'_{i}$ is any fixed subspace of $V_{i} \in \mathcal{V}$ with  $dim(V'_{i}) = s$
if $dim(V_{i}) \geq s$, and $V'_{i} = V_{i} \in \mathcal{V}$ if $dim(V_{i}) < s$. Then we have $|\mathcal{V}| = |\mathcal{V}'|$.
Since $\mathcal{V}$ is $\mathcal{L}$-intersecting Sperner, $\mathcal{V}'$ is $\mathcal{L}$-intersecting Sperner. It follows from Theorem \ref{thm3.4} that
\[m = |\mathcal{V}| =|\mathcal{V}'| \leq \qbinom{n}{s}.\]

\noindent {\bf Case 2.}
$\mathcal{L} \neq \{0, 1, \dots, s-1\}$.

Let $\mathcal{V}' = \{V_{j} \in \mathcal{V} \mid dim(V_{j}) \leq s\}$ and $\mathcal{L}' = \mathcal{L} \cap \{0, 1, \dots, s - 1\}$. Then $|\mathcal{L}'| \leq s - 1$.
Clearly, $\mathcal{V}'$ is $\mathcal{L}'$-intersecting. By the induction hypothesis, we have
\[|\mathcal{V}'| \leq \qbinom{n}{s-1}.\]
Denote by $\qbinom{\mathcal{W}}{1}$ the set of all $1$-dimensional subspaces of the vector space $\mathcal{W}$. Then
\[ \left|\qbinom{\mathcal{W}}{1}\right| = \qbinom{n}{1}.\]
For each $1$-dimensional subspace $R$  of $\mathcal{W}$, let
\[\mathcal{V}_{R} = \{V_{j} \in \mathcal{V} \mid R \subseteq V_{j}\}.\]
Then $\mathcal{V}_{R}$ is $\mathcal{L}^{*}$-intersecting, where $\mathcal{L}^{*} = \mathcal{L}\setminus \{0\} = \{l_{2}, \dots, l_{s}\}$.
Set
\[\mathcal{V}'_{R} = \{V_{j}/R \mid V_{j} \in \mathcal{V}_{R}\},\]
where $V_{j}/R$ is the quotient of two subspaces $V_{j}$ and $R$.
Then it is clear that $\mathcal{V}'_{R}$ is an $\mathcal{L}^{**}$-intersecting family of subspaces of the $(n-1)$-dimensional vector space
$\mathcal{W}/R$, where $\mathcal{L}^{**} = \{l_{2} - 1, \dots, l_{s} - 1\}$.
Since $\mathcal{V}$ is Sperner,  $\mathcal{V}'_{R}$ is Sperner.
By the induction hypothesis, we have that, for each $1$-dimensional subspace $R$  of $\mathcal{W}$,
\[|\mathcal{V}_{R}| = |\mathcal{V}'_{R}|  \leq \qbinom{n-1}{s-1}.\]
Note that for each $V_{j} \in \mathcal{V} \setminus \mathcal{V}'$, $dim(V_{j}) \geq s + 1$. It follows that
\[(m - |\mathcal{V}'|)\qbinom{s+1}{1} \leq \sum_{R \in \qbinom{\mathcal{W}}{1}}|\mathcal{V}_{R}| \leq \qbinom{n}{1} \qbinom{n-1}{s-1}.\]
Recall that for $n \geq 2s + 1$,
\[\qbinom{n}{s-1} = q^{s - 1}\qbinom{n-1}{s-1} + \qbinom{n-1}{s-2}
\leq (q^{s - 1} + 1)\qbinom{n-1}{s-1}.\]
Since $q \geq 2$ and $n \geq 2s + 1$, it follows that
\begin{align*}
m &= |\mathcal{V}| = |\mathcal{V}'| + (m - |\mathcal{V}'|)   \\
&\leq \qbinom{n}{s-1} + \frac{\qbinom{n}{1}}{\qbinom{s+1}{1}}\qbinom{n-1}{s-1} \\
&\leq (q^{s - 1} + 1)\qbinom{n-1}{s-1}
 + \frac{q^{n} - 1}{q^{s+1} - 1} \qbinom{n-1}{s-1}  \\
& = \left((q^{s - 1} + 1) + \frac{q^{n} - 1}{q^{s+1} - 1}\right)\qbinom{n-1}{s-1}  \\
& < \left((q^{s - 1} + 1) + \frac{q^{n} - 1}{q^{s+1} - q}\right)\qbinom{n-1}{s-1}  \\
& = \left(\frac{(q^{s - 1} + 1)q(q^{s} - 1)}{q(q^{s} - 1)} + \frac{q^{n} - 1}{q(q^{s} - 1)}\right)\qbinom{n-1}{s-1} \\
& < \left(\frac{q^{2s+1}-1}{q(q^{s} - 1)} + \frac{q^{n} - 1}{q(q^{s} - 1)}\right)\qbinom{n-1}{s-1}  \\
& \leq \left(\frac{q^{n}-1}{q(q^{s} - 1)} + \frac{q^{n} - 1}{q(q^{s} - 1)}\right)\qbinom{n-1}{s-1}  \\
&\leq \frac{q^{n} - 1}{q^{s} - 1}\qbinom{n-1}{s-1}
= \qbinom{n}{s}.
\end{align*}
$\hfill \Box$

\section{Concluding Remarks}
In section 2, we verified Snevily's conjecture for $n \geq {{k^{2}} \choose {l_{1}+1}}s + l_{1}$ with $k$ being the maximum subset size (Theorem 1.9) and derived a $t$-wise extension (Theorem 1.10). We then provided a vector space analogue to Snevily's conjecture (Theorems \ref{thm1.15} and \ref{thm1.16}) in section 3.
An open problem is to determine whether Theorems \ref{thm1.9}, \ref{thm1.15}, and \ref{thm1.16}
hold for small values of $n$ between $2s+1$ and the bounds in the corresponding theorems.

By Lemma \ref{lem3.1} and Theorem \ref{thm3.3}, one has the next well-known Sperner theorem for families of subspaces of an $n$-dimensional vector space $\mathcal{W}$ over finite field $\mathbb{F}_{q}$, which is a vector space analogue of the classical Sperner theorem.

\begin{theorem}\label{thm4.1} (q-analogue Sperner Theorem).  Assume that $\mathcal{V}$ is a Sperner family
of subspaces of an $n$-dimensional vector space over a finite field $\mathbb{F}_{q}$. Then
\[|\mathcal{V}| \leq \qbinom{n}{\lfloor \frac{n}{2}\rfloor}.\]
\end{theorem}

In addition to strengthening Theorems  \ref{thm1.12} -- \ref{thm1.14}, Theorem \ref{thm1.16} also shows that the upper bound in the q-analogue Sperner Theorem (Theorem 4.1) can be improved under $\mathcal{L}$-intersecting condition when $n \geq 2s + 1$.

\section*{Declaration of Competing Interest}
The authors declare that they have no conflicts of interest to this work.

\section*{Acknowledgement}
The research is supported by the National Natural Science Foundation of China (71973103, 11861019, 12271527, 12071484).
Guizhou Talent Development Project in Science and Technology (KY[2018]046), Natural Science Foundation of Guizhou ([2019]1047, [2020]1Z001,  [2021]5609).

\end{document}